\documentclass{amsart}

\usepackage{amssymb}
\usepackage{amsmath}
\usepackage{amsthm}
\usepackage{amsbsy}

\theoremstyle{plain}
\newtheorem{theorem}{Theorem}

\newtheorem{lemma}[theorem]{Lemma}

\newtheorem{proposition}[theorem]{Proposition}

\theoremstyle{definition}

\theoremstyle{remark}

\newcommand{\del}{\partial}

\newcommand{\sm}{\setminus}

\newcommand{\co}{\colon\thinspace}

\begin{document}

\title{A proof of a homeomorphism theorem of Waldhausen }

\author{Siddhartha Gadgil}
\address{       Department of mathematics\\
                SUNY at Stony Brook\\
                Stony Brook, NY 11794}
\email{gadgil@math.sunysb.edu}

\author{Gadde A. Swarup} 
\address{Department of Mathematics and Statistics\\
	 University of Melbourne\\
	 Parkville, Victoria\\
	 Australia} 

\email{gadde@ms.unimelb.edu.au}

\begin{abstract}
We give a short proof of Waldhausen's homeomorphism theorem for
orientable Haken manifolds.
\end{abstract}

\maketitle

We give here a short proof of Waldhausen's homeomorphism
theorem~\cite{Wa} for orientable Haken manifolds.

Recall that a compact irreducible 3-manifold is Haken (i.e., it has a
hierarchy) if and only if it is sufficiently large (i.e., it admits an
incompressible surface).  This follows from the Haken-Kneser argument
which shows that there is a bound on the number closed disjoint
non-parallel incompressible surfaces in a compact irreducible
3-manifold. A readable proof of this is given in Hempel's
book~\cite{He}.

We now state Waldhausen's homeomorphism theorem. The proof that
follows is a variant of the standard arguments.  In the standard
proofs, some of the special cases of Waldhausen's theorem are proved
during the course of the proof. We avoid this by appealing to a result
of Hopf and Kneser that algebraic and geometric degrees are equal 
(see~\cite{Ep} and~\cite{Sk}). 

We will always assume that we are in the orientable case and that all
our submanifolds are two sided. The argument also works in the case of
non-orientable Haken manifolds.

\begin{theorem}
Suppose $(M,\del M)$ and $(N,\del N)$ are compact irreducible
$3$-manifolds with (possibly empty) boundary and $N$ is sufficiently
large. Then any homotopy equivalence $f\co (M,\del M)\to (N,\del N)$,
whose restriction to $\del M$ is a homeomophism, is homotopic to a
homeomorphism.
\end{theorem}

We first consider the case when $\del N$ is non-empty. Then there
exists an incompressible, $\del$-incompressible surface $S\subset N$
such that $\del S\neq\phi$. Let $f^{-1}(S)=\cup_{i=1}^n F_i$, where
each $F_i$ is connected. The following lemma is by now standard and
stems from the work of Whitehead and Stallings.

\begin{lemma}
After a homotopy of $f$, we may assume that each $F_i$ is essential.
\end{lemma}

Henceforth assume that each $F_i$ is essential. Further, we assume
that $f$ is a homeomorphism in a neighbourhood of $\del M$.

\begin{proposition}
$f^{-1}(S)$ is connected and $f\co f^{-1}(S)\to S$ is degree-one.
\end{proposition}

\begin{proof}
Let $f_k\co F_k\to S $, $1\leq k\leq m$ be the restrictions of $f$. We
shall use the fact that the degree of a map can be computed locally.

Thus, let $\del S=\gamma_1\cup\dots\cup\gamma_m$. Then, as $f$ is a
homeomorphism in a neighbourhood of the boundary, for $p$ in a
neighbourhood of any $\gamma_i$, the inverse image $f^{-1}(p)$
consists of a single point. It follows that

$$deg(f_k)=\begin{cases}
	\pm 1 & \text{if $f^{-1}(\gamma_i)\subset F_k$}\\
	0 & \text{otherwise}
	\end{cases}$$

It follows that exactly one $f_k$, say $f_1$ has degree $\pm 1$ and
the others have degree $0$. Further, it follows that if $k\neq 1$
$\del F_k=\phi$. But this implies that $F_k$ is a closed surface. As
$S$ has a non-empty boundary, $f_k$ is not an injection on $\pi_1$. As
$f$ is a homotopy equivalence and $S$ is essential, $F_k$ cannot be
essential, a contradiction. Thus $f^{-1}(S)=F_1$ is connected and
$f_1$ has degree $\pm 1$.
\end{proof}

Next, we shall show that the same result can be achieved in the case
when $\del M=\phi$. We use a theorem of Hopf and Knesser, which says
that the algebraic and geometric degrees are equal (for a modern proof
see~\cite{Sk}). As $f$ is degree-one, given a ball $B\subset N$, after
a homotopy of $f$, the restriction of $f$ to $B'=f^{-1}(B)$ is a
homeomorphism.

Let $S$ be an essential surface in $N$. By an isotopy of $S$, we can
ensure that $D=S\cap B$ is a properly embedded disc in $S$. Further,
we may apply the lemma to ensure that $f^{-1}(S)$ is essential without
altering this. Using the notation as in the previous case, we once
more have the following.

\begin{proposition}
$f^{-1}(S)$ is connected and $f\co f^{-1}(S)\to S$ is degree-one.
\end{proposition}
\begin{proof}
Again, using the fact that degree can be computed locally, and picking
$p\in D$, we get that, after possibly re-ordering components, $f_1$
has degree $\pm 1$ and all other $f_k$ have degree $0$. Further, if
$k\neq 1$, the image of $f_k$ lies in $S\sm D$, and hence $f_k$ cannot
inject in $\pi_1$, giving a contradiction as before.
\end{proof}

Thus, in either case $F=F_1=f^{-1}(S)$ is connected and $f_{1*}$ is an
injection on fundamental groups. As it has degree one, it is also a
surjection. Further, its restriction to $\del F_1$ is a
homeomorphism. It follows that $f_1$ is homotopic to a homeomorphism.

Thus we can cut $M$ and $N$ and use finiteness of hierarchies to
complete the proof.

\end{document}